\documentstyle{article}
\renewcommand{\baselinestretch}{1.3}

\newtheorem {th}{Theorem}
\newtheorem {lem}[th]{Lemma}
\newtheorem {pr}[th]{Proposition}
\newtheorem {cor}[th]{Corollary}
\newtheorem{defn}{Definition}

\def\Cox{\hfill \Box \vspace{2ex}}
\def\deq{\, {\stackrel {def} {=}}}

\def\ul{\underline}

\def\ee{\epsilon}

\def\iid{\mbox{i.i.d.}}

\def\E{{\bf{E}}}
\def\P{{\bf{P}}}

\def\R{{\bf{R}}}
\def\Z{{\bf{Z}}}
\def\S{{\cal{S}}}

\def\|{\, | \, }

\def\Cap{\mbox{Cap}}

\begin{document}
 
\begin{center}
{\large \bf CRITICAL RANDOM WALK IN RANDOM ENVIRONMENT ON TREES
OF EXPONENTIAL GROWTH} \\
\end{center}
\vspace{5ex}
\begin{center}
ROBIN PEMANTLE \footnote{Research supported in part by a
National Science Foundation postdoctoral fellowship
and be NSF grant \# DMS9103738} \footnote{Department of
Mathematics, University of Wisconsin-Madison, Van Vleck Hall, 480 Lincoln
Drive, Madison, WI 53706}
\end{center}

\vspace{4ex}

{\bf ABSTRACT:} \break
This paper studies the behavior of RWRE on trees in the critical case 
left open in previous work.  For trees of exponential growth, a random
perturbation of the transition probabilities can change a transient
random walk into a recurrent one.  This is the opposite of what occurs
on trees of sub-exponential growth.  
\vspace{5ex}

\section{Introduction}

This paper is concerned with the problem of determining
whether a random walk in a random environment (RWRE) on an infinite,
exponentially growing tree
is transient or recurrent.  The problem was first studied in
\cite{Pe} as a way of analyzing another process called 
Reinforced Random Walk, and then in \cite{LP} where a more
complete solution was obtained.  It was shown there that the RWRE
is transient when the size of the tree, as measured by the log of the
{\em branching number}, is greater than the {\em backward
push} of the random environment, and recurrent when the
log of the branching number is smaller than the backward push. The 
case of  equality was left open.  For trees of sub-exponential growth,
this critical case was almost completely settled in \cite{PP}.
The present paper is a companion to \cite{PP} in that it attempts
to settle the critical case for exponentially growing trees.
The results here are less definitive than in the sub-exponential
case, in that the sufficient conditions for transience and for
recurrence are in terms of capacity and growth respectively;
these conditions are not quite complementary,
leaving open a critical-within-critical case.  A technical assumption
on the random environment is also required;
examples show that this assumption is often satisfied.  On the
positive side, it is proved here that a phase boundary occurs
in an unusual place, namely when the growth rate of the tree, 
$\exp (\beta n + o(n))$, has the $o(n)$ term equal to a  constant
times $n^{1/3}$.  Here follows a precise description of the problem.

Let $\Gamma$ be any infinite, locally finite tree with no leaves
(vertices of degree one).  Designate a vertex $\rho$ of
$\Gamma$ as its root.  For any vertex $\sigma \neq \rho$, denote by 
$\sigma'$ the unique neighbor of $\sigma$ closer to $\rho$ ($\sigma'$ is
also called the parent of $\sigma$).  An environment for random walk on
a fixed tree, $\Gamma$, is a choice of transition probabilities 
$q(\sigma, \tau)$ on the vertices of $\Gamma$ with $q(\sigma , \tau) > 0$
if and only if $\sigma$ and $\tau$ are neighbors.  When these transition
probabilities are taken as random variables, the resulting mixture of
Markov chains is called Random Walk in Random Environment (RWRE).
Following \cite{LP} and the references therein, random environments 
studied in this paper satisfy the homogeneity condition 
\begin{equation} \label{eq1}
\mbox{The variables } X (\sigma) = \ln \left ( {q(\sigma' , \sigma) 
   \over q(\sigma' , \sigma'')} \right ) \mbox{ are i.i.d.\ for } |\sigma|
   \geq 2 ,
\end{equation}
where $|\sigma|$ denotes the distance from $\sigma$ to $\rho$.  Here, and
throughout, let $X$ denote a random variable with this common distribution.  

Before stating the main result, a few definitions and notations are required.  
Write $\sigma \leq \tau$ if $\sigma$ is on the path connecting $\rho$ and
$\tau$; in this paper, the term ``path'' always refers to a 
path without self-intersection.
Write $\sigma \wedge \tau$ for the greatest lower bound
of $\sigma$ and $\tau$; pictorially, this is where the paths from
$\rho$ to $\sigma$ and $\tau$ diverge.  Let $\partial \Gamma$, called
the {\em boundary} of $\Gamma$, denote the set of infinite 
paths beginning at $\rho$.  Let $\Gamma_n$ denote the set
$\{ \sigma : |\sigma| = n \}$ of vertices at level $n$ of $\Gamma$.
Define the {\em backward push} of the random environment, denoted
$\beta$ or $\beta (X)$, by
$$ \beta (X) = - \ln \min_{0 \leq \lambda \leq 1} \E e^{\lambda X} .$$
The size of an infinite tree is best discussed in terms of capacity.
\begin{defn}
Let $\phi : \Z^+ \rightarrow \R^+$ be a nonincreasing function.  Define
the $\phi$-energy of a probability measure $\mu$ on the boundary of
$\Gamma$ to be
$$I_\phi (\mu) = \int_{\partial \Gamma} \int_{\partial \Gamma}
   \phi (|\xi \wedge \eta|)^{-1} \, d\mu (\xi) \, d\mu (\eta) .$$
Define the capacity of $\Gamma$ in gauge $\phi$ by
$$\Cap_\phi (\Gamma) = \left [ \inf_\mu I_\phi (\mu) \right ] ,$$
where the infimum is over all probability measures on $\partial \Gamma$.,
and $\Cap_\phi (\Gamma) \neq 0$ if and only if there is some measure
of finite energy.
\end{defn}
Say that $\Gamma$ is {\em spherically symmetric} if there is a
{\em growth function} $f: \Z^+ \rightarrow \Z^+$ such that every
vertex $\sigma \neq \rho$ has $1 + f(|\sigma|)$ neighbors; in other
words, the degree of a vertex depends only on its distance from the
root.  A spherically symmetric tree $\Gamma$ 
has positive capacity in gauge $\phi$ if and only if
$$\sum \phi (n) |\Gamma_n|^{-1} < \infty .$$
Thus positive capacity in gauge $\phi (n) = e^{-kn}$ implies
liminf exponential growth rate of at least $k$.  In particular,
the supremum of those $k$ for which $\Gamma$ has positive capacity
in gauge $\phi (n) = e^{-kn}$ is the Hausdorff dimension, $dim (\Gamma)$;
in the terminology of~\cite{Ly1} and~\cite{LP}, $dim(\Gamma)$ is
the log of the branching number.

The main result of \cite{LP} is that RWRE on $\Gamma$ is
a.s.\ transient if $dim (\Gamma) > \beta (X)$,
and a.s.\ recurrent if $dim(\Gamma) < \beta (X)$.  The use of
gauges more general than $e^{-kn}$ allows for finer distinctions
of size to be made within the class of trees of the same
dimension.  In particular, when $\beta = dim (\Gamma) = 0$,
it is shown in \cite{PP} that positive capacity in gauge
$n^{1/2}$ is sufficient and almost necessary for transience
of RWRE.  Since the case $\beta = 0$ is in some sense a mean-zero 
perturbation of the deterministic environment of
a simple random walk ($X \equiv 0$), and simple random walk
is transient if and only if $\Gamma$ has positive capacity in
gauge $n^{-1}$, this shows that the perturbation makes the
walk more transient.  By contrast, the main result of this
paper is as follows. 

\begin{defn} \label{deftopheavy}
Say that a real random variable 
$X$ is {\em top-heavy} if the infimum of $\E e^{\lambda X}$
over $\lambda \in [0,1]$ is achieved at some $\lambda_0 \in (0,1)$
and $\E e^{\gamma X} < \infty$ for some $\gamma > \lambda_0$.
\end{defn}

\begin{th} \label{th2.3}
Consider RWRE on a tree $\Gamma$ with $\beta (X) = dim(\Gamma) > 0$.
If $X = -\beta$ with probability one, then RWRE is transient if and only
if $\Gamma$ has positive capacity in gauge $\phi (n) = e^{-n \beta}$.
On the other hand if $X$ is nondeterministic, top-heavy, 
and is either a lattice distribution or has an absolutely continuous
component with density bounded above and bounded away from zero in a 
neighborhood of zero, then
\begin{quote}
$(i)$  there exists $c_1 (X)$ for which the growth bound
$$|\Gamma_n| \leq e^{\beta n + c_1 n^{1/3}} \mbox{ for all } n$$ 
implies that RWRE is recurrent;

$(ii)$  there exists $c_2 (X)$ such that if $\Gamma$ has positive
capacity in gauge $\phi (n) = e^{- n \beta - c_2 n^{1/3}}$ then
RWRE is transient.  
\end{quote}
\end{th}

\noindent{{\bf Remark:}} The requirement that $X$ be 
top-heavy is enigmatic, but
not overly restrictive.  For example, it is satisfied by 
normal random variables with mean $-c$ and variance $V$ 
whenever $c < 2V$.  In the case where $X$ takes only the 
values $\pm 1$, with $\P (X = 1) = p < 1/2$, it is 
top-heavy if and only if $(1-p)/p < e^2$.  

The remainder of this section outlines the the proof of this theorem
and serves as a guide to the remaining sections.  
Theorem~\ref{th2.3} is proved in the following three steps.  First, 
in Section~2,
a correspondence connection between random walks and electrical networks 
\cite{DS} reduces the problem to one of determining whether a random 
electrical network is transient or recurrent almost surely. After
this reduction, the technical condition of top-heaviness comes in: 
top-heaviness implies that finite resistance, when achieved, will
be due to a single random infinite path of finite resistance; 
searching for a single path with this property is easier than 
searching for some large collection of paths with a weaker property.
Next, large deviation estimates are needed for the probability of
an unusually small resistance along a fixed path of length $n$
(Lemma~\ref{lem tube estimate} and Corollary~\ref{cor1}).
These are applied via a simple first-moment calculation to
obtain Lemma~\ref{lem small prob}, which is
an upper bound tending to zero on the probability that any of the
$|\Gamma_n|$ chains of resistances of length $n$ stays small.
Incidentally, this is where the ``extra'' factor of $e^{cn^{1/3}}$
comes in.  Estimates with this same factor have been obtained
by Kesten with much greater accuracy for branching Brownian
motion~\cite{Ke}.  Ours are a discrete analogue of Kesten's in
the sense that continuous-time branching has been replaced
by $\Gamma$-indexed branching; this analogy is explained
more fully in~\cite{BP} and~\cite{PP}.
Part~$(i)$ of Theorem~\ref{th2.3} follows from this upper bound 
by computing the expected truncated conductance.  The last step, 
which is needed only for the proof of part~$(ii)$, is a 
second-moment technique (Lemma~\ref{lem2mm}) developed 
in \cite{Ly2} and \cite{PP} for proving the almost sure
existence of an infinite path of finite resistance based on
the two-dimensional marginals for finite paths (i.e.\ the
probabilities for two paths of length $n$ both to have 
large conductance if the paths share the first $k$ resistors). 

\section{Reductions}

Begin with the reduction of the recurrence/transience problem to
an electrical problem.  As is well known, transience of a reversible
Markov chain is equivalent to finite resistance of the 
associated resistor network on the same graph, where the transition 
probabilities
from any vertex are proportional to the conductances (reciprocal
resistances) of the edges incident to that vertex; see for example
\cite{DS}.  For a random environment satisfying~(\ref{eq1}),
the resistances in the associated random electrical network
are easily seen to be given by
$$\mbox{Resistance along } \overline{ \sigma' \sigma} = e^{-S(\sigma)} $$
where 
$$S(\sigma) = \prod_{\rho < \tau \leq \sigma} X(\sigma) .$$
Here, the values of $X(\sigma)$ for $|\sigma| \leq 1$ are 
assigned to make this relation hold for all $\sigma \neq \rho$, 
while the values for $|\sigma| \geq 2$ are i.i.d.\ by~(\ref{eq1}).
Since finiteness of the total resistance is not affected by changing 
finitely many resistances, we alter the $X(\sigma)$ for 
$|\sigma| \leq 1$ so that the entire collection 
is i.i.d.  

A sufficient condition for transience is the existence of an 
infinite path $\rho, \sigma_1 , \sigma_2 , \ldots$ along which
$\sum e^{- S(\sigma)} < \infty$.  Conversely, let
$$U (\sigma) = \min_{\rho < \tau \leq \sigma} e^{S(\sigma)} .$$
A useful sufficient condition for recurrence is given by
the following lemma.

\begin{lem} \label{bottleneck}
Let $\Gamma$ be any tree with conductances $C(\sigma)$, and
let $\Pi$ be any cutset, i.e.\ any minimal set among those
intersecting every infinite path from $\rho$.  Then
the conductance from $\rho$ to $\Pi$ is at most 
$$ \sum_{\sigma \in \Pi} U(\sigma) .  $$
Consequently, if the conductances are random with  
$\sum_{|\sigma| = n} U(\sigma) \rightarrow 0$
in probability, then the random walk is recurrent with
probability one.
\end{lem}

\noindent{{\bf Proof:}}  For each $\sigma \in \Pi$, let $\gamma (\sigma)$
be the sequence of conductances on the path from $\rho$ to $\sigma$, 
and let $\Gamma'$ be a tree consisting of disjoint paths for each $\sigma 
\in \Pi$, each path having conductances $\gamma (\sigma)$.  $\Gamma$
is a contraction of $\Gamma'$, so by Rayleigh's monotonicity law, the
conductance to $\Pi$ in $\Gamma$ is less than or equal to the conductance
of $\Gamma'$, which is the sum over $\sigma \in \Pi$ of conductances
bounded above by $U(\sigma)$.  Putting $\Pi = \Gamma_n$ shows that
the conductance from $\rho$ to infinity is bounded above by
$\liminf_n \sum_{|\sigma| = n} U(\sigma)$, proving the lemma.   $\Cox$

In the next section estimates will be given on $\P (\sigma \in W)$
and $\P (\sigma , \tau \in W)$, where for fixed constants $c$ and $L$,
$W$ is the set of vertices $\sigma$ such that for every 
$\tau \leq \sigma$ with $|\tau| > L$, 
$$1/10 \leq S(\tau) / c |\tau|^{1/3} \leq 1.$$
These estimates are then plugged into the
following result of Lyons (see \cite[Theorem~4.1]{PP}).

\begin{lem} \label{lem2mm}
Let $\Gamma$ be an infinite, locally finite tree without leaves,
let $X(\sigma)$ be i.i.d.\ random variables indexed by the vertices
of $\Gamma$ and let $B_n$ be a subset of $\R^n$ for each $n$.
Let $W$ be the set of vertices $\sigma \in \Gamma$
such that for every $\tau \leq \sigma$, the sequence $(X(\rho) , 
\ldots , X(\tau))$ along the path from $\rho$ to $\tau$ is in
the set $B_{|\tau|}$.  Suppose there is a positive, nonincreasing
function $g : \Z^+ \rightarrow \R$ such that for any two vertices
$\sigma, \tau \in \Gamma_n$ with $|\sigma \wedge \tau| = k$,
\begin{equation} \label{eq2mm}
\P (\sigma , \tau \in W) \leq {\P (\sigma \in W)^2 \over g(k)} .
\end{equation}
Then the probability of $W$ containing an infinite path is
at least $\Cap_g(\Gamma)$.    $\Cox$
\end{lem}

\noindent{{\bf Remark:}} Of course any infinite path in $W$ has 
$\sum e^{-S(\sigma)} < \infty$, implying transience.  The
reason that one looks for a path along which $S(\tau)$ 
is bounded above as well as below is so as to be able to apply 
this lemma, which is really a
jazzed up second moment bound.  It is important to find
a random subset of $\Gamma_n$ whose cardinality has a second
moment not too much larger than the square of
its first moment; then the LHS of~(\ref{eq2mm}) will
not be too large.    

\section{Large deviation estimates}

The first series of estimates concern the probabilities
$\P (\sigma \in W)$ of the previous lemma.  The routine proofs of
the first two propositions in the series are omitted.

\begin{pr} \label{brownian}
Let $\P_x$ be the law of a standard one-dimensional Brownian motion
started at $x$.  Let $-1 <  c_1 < c_2 < 1$ and $-1 < c_3 < c_4 < 1$
be real constants.  Then there exist positive constants $K_1$ and $K_2$ 
such that for any $L > 0$, the following two inequalities hold.
\begin{eqnarray*}
\sup_x \P_x (|B_t| \leq 1 \mbox{ for all } t \leq L) & \leq & K_1 e^{-{\pi^2
   \over 8} L} \hspace{.2in} ; \\[2ex] 
   \inf_{c_1 \leq x \leq c_2} \P_x (|B_t| \leq 1 \mbox{ for all } t \leq L
   \mbox{ and } c_3 \leq |B_L| \leq c_4) & \geq & K_2 e^{-{\pi^2 \over 8} L} .
\end{eqnarray*}
\end{pr}

\begin{pr} \label{rw}
Let $\{ X_n \}$ be a sequence of i.i.d.\ random variables with mean zero and 
variance $V < \infty$.  Let $S_n = \sum_{i=1}^n X_i$ and for $0 < t <1$,
interpolate polygonally by letting $S_{n+t} = S_n + tX_{n+1}$.   Then
\begin{eqnarray} 
&& \limsup_{L \rightarrow \infty} \limsup_{n \rightarrow \infty}   
   \sup_{|x| \leq 1} L^{-1} \, \ln \left ( \P ( |(Vn)^{-1/2} S_{nt} + x| 
   \leq 1 \mbox{ for all } t \leq L) \right ) \nonumber \\[2ex]
& & \leq - {\pi^2 \over 8} \label{upper rw}.
\end{eqnarray} 
If $-1 <  c_1 < c_2 < 1$ and $-1 < c_3 < c_4 < 1$ are real constants, then
\begin{eqnarray} 
&& \liminf_{L \rightarrow \infty} \liminf_{n \rightarrow \infty}   
   \inf_{c_1 \leq x \leq c_2} L^{-1} \,  \ln \P [ A(L,n,x)] \nonumber \\[2ex]
& \geq & - {\pi^2 \over 8} , \label{lower rw}
\end{eqnarray} 
where $A(L,n,x)$ is the event
$$ \left \{ |(Vn)^{-1/2} S_{nt} + x| \leq 1 \mbox{ for all } t \leq L,
   \mbox{ and } c_3 \leq S_{nL} + x \leq c_4 \right \} .$$
\end{pr}

\begin{lem} \label{lem tube estimate}
Suppose $f , g : \Z^+ \rightarrow \R$ satisfy $f > g$ and 
$\lim_{t \rightarrow \infty}
f(t) - g(t) = \infty$, and assume the following flatness hypothesis:
$$ \sup_L \lim_{t \rightarrow \infty} \; {\sup_{0 \leq s \leq L(f(t) 
   - g(t))^2}  \max (|f(t+s) - f(t)| , |g(t+s) - g(t)|) \over f(t) - g(t)}
   \; = \, 0 .$$
Let $S_n = \sum_{i=1}^n X_i$ be a random walk with $\E X_1 = 0$ and
$\E X_1^2 = V < \infty$ and suppose that for each $n$, $\P (g(k) < S_k < f(k)
\mbox{ for all } k \leq n) > 0$.  Then
\begin{eqnarray} 
&&  \lim_{n \rightarrow \infty} \left ( \sum_{k=1}^n (f(k) - g(k))^{-2} 
   \right )^{-1} \ln \left ( \P (g(k) < S_k < f(k)
   \mbox{ for all } k \leq n) \right ) \nonumber \\[2ex]
& & = \; - {\pi^2 \over 8} V. \label{eqth1}
\end{eqnarray}
\end{lem}

\begin{cor} \label{cor1}
Let $Y_n$ be i.i.d.\ mean zero random variables with partial sums 
$T_n = \sum_{i=1}^n Y_i$ and suppose that $\E Y_1^2 = V < \infty$.
Then for any real $c_1< c_2$,  
$$\lim_{n \rightarrow \infty} \; {\ln \P \left ( c_1 k^{1/3} \leq T_k
    \leq c_2 k^{1/3} \mbox{ for all } k \leq n \right ) \over n^{1/3}}
   = -{\pi^2 \over 8} {3V \over (c_2 - c_1)^2} $$
unless the probability is eventually zero.  
\end{cor}

\noindent{{\bf Proofs:}}  The first proposition is a standard eigenvalue
estimate  -- see for example \cite[p.\ 259]{Du}.  The second
follows from the first and the invariance principle by piecing 
together segments of path of length $L$; the two-page proof
is entirely routine and is omitted.  

To establish the lemma, prove first that the limsup is at most 
$-{\pi^2 \over 8} V$.  Extend $f$, $g$ and $S$
to $\R^+$ by linear interpolation with $f(0) = g(0) = S_0 = 0$; for any
integer $n$, the event $\{ g(k) \leq S_k \leq f(k) : k = 1 , \ldots , n \}$
is the same when $k$ is replaced by a real parameter $t$ running from 
0 to $n$.  Fix $L > 0$. For any positive integer $m_0$ inductively 
define $m_{k+1} = \lceil m_k + (1/4)V^{-1}L(f(m_k)- g(m_k))^2 \rceil$.  
Let $M(n)$ denote $\max \{ k : m_k < n \}$.
\begin{quote}
\ul{Claim}:  $\lim_{n \rightarrow \infty} M(n)^{-1} \sum_{k=1}^n
4 VL^{-1} (f(k) - g(k))^{-2} = 1$.
\end{quote}

\noindent{To} prove this, first choose $\ee \in (0,1/8)$ and $j$ large enough 
so that $t \geq m_j$ implies
\begin{equation} \label{slow var}
\sup_s { \max (|f(t+s) - f(t)| , |g(t+s)-g(t)|) \over f(t) - g(t)} < \ee ,
\end{equation}
where the supremum is over $s$ for which $0 \leq s \leq (1/4) 
V^{-1}L (f(t) - g(t))^2$.  Summing the identity 
$\sum_{i = m_k}^{m_{k+1}} (m_{k+1} - 
m_k)^{-1} = 1$ over $j \leq k \leq M(n) + 1$ yields
\begin{eqnarray*} 
M(n) & \leq & j +1 + \sum_{i = m_0}^n 4VL^{-1} (f(m_{M(i)}) - 
   g(m_{M(i)}))^{-2} \\[2ex]
& \leq &  j + 1 + \sum_{i = m_0}^n 4VL^{-1} (f(i) - g(i))^{-2} 
   (1-2\ee)^{-2} .
\end{eqnarray*} 
Similarly,
$$ M(n)  \geq  j + \sum_{i = m_0}^n 4VL^{-1} (f(i) - g(i))^{-2}
    (1+2\ee)^{-2} .$$
Letting $n \rightarrow \infty$ proves that the limsup and liminf of
$$M(n)^{-1}  \sum_{k=1}^n 4VL^{-1} (f(k) - g(k))^{-2}$$
are between $(1+2\ee)^{-1}$ and $(1-2\ee)^{-1}$.  Letting 
$\ee \rightarrow 0$ proves the claim.  

Continuing the proof of the lemma, pick $\ee > 0$ and $N$ large enough
so that for our (still fixed) value of $L$ and $n \geq N$, the
sup in~(\ref{upper rw}) is at most $-{\pi^2 \over 8} + \ee$.  Pick $m_0$ large 
enough so that $t \geq m_0$ implies firstly~(\ref{slow var}) and secondly
$f(t) - g(t) > (4VN)^{1/2}$.  Let $y_k = (f(m_k) + 
g(m_k))/2$ and $d_k = (f(k) - g(k))/2$, so $m_{k+1}
- m_k = \lceil LV^{-1} d_k^2 \rceil$.  Let $R_k$ be the rectangle 
$$\{ (x,y) : m_k \leq x \leq m_{k+1} \mbox{ and } y_k - (1+ 2\ee) d_k 
\leq y \leq y_k + (1+ 2\ee) d_k \} .$$  
Then the union of the $R_k$ covers the graphs of $f$ and $g$ after
$m_0$, i.e.\ the region $ \{ (x,y) : x \geq m_0 , g(x) \leq y \leq f(x) \}$. 
Write $S_i'$ for a random walk independent of $S_i$ and identically 
distributed.  Using the Markov property of $S_i$, the stationarity of
its increments, and then rescaling each $R_k$ to have height one, gives
\begin{eqnarray*}
&& \P (f(k) \leq S_k \leq g(k) : 1 \leq k \leq n) \\[2ex]
& \leq & \P ((i,S_i) \in \bigcup R_k : m_0 \leq i \leq n) \\[2ex]
& \leq & \prod_{k=0}^{M(n)-1} \P \left ( (i,S_i) \in R_{k-1} : m_k \leq 
   i \leq m_{k+1} \| S_{m_k} \right ) \\[2ex]
& = & \prod_{k=1}^{M(n)-1} \P (|S_{m_k} + S_i' - y_k| \leq (1+2\ee)d_k    
   : i = 1 , \ldots , m_{k+1} - m_k) \\[2ex]
& = & \prod_{k=1}^{M(n)-1} \P \left ( |(1+2\ee)^{-1} d_k^{-1} 
   [(S_{m_k} - y_k) + S_{V^{-1} (1+2\ee)^2 d_k^2 t}'] | \leq 1 \right. \\
&& \left. \hspace{.5in} \mbox{ for all }
   0 \leq t \leq L (1+2\ee)^{-2} \right ) \\[2ex]
& \leq & \sup_{w \geq N} \, \sup_{0 \leq x \leq 1} \left [ \P (
    {|x + S_{wt}| \over (Vw)^{1/2}} \leq 1 : 0 \leq t \leq L (1+2\ee)^{-2}) 
    \right ]^{M(n)-1}
\end{eqnarray*}
where the rescaling factor $w = (1+2\ee)^2 V^{-1} d_k^2$
is at least $N$ by the previous choice
of $m_0$.   Taking the log, dividing by $L(1+2\ee)^{-2}$ and
applying~(\ref{upper rw}) gives
$$L^{-1}(1+2\ee)^2  \ln \P (f(k) \leq S_k \leq g(k) : k \leq n) 
   \leq (M(n) - 1) (-{\pi^2 \over 8} + \ee)$$
by choice of $N$.  Plugging in the asymptotic value of $M(n)$ from the 
claim above gives that for sufficiently large $n$,
\begin{eqnarray*}
&& \ln \P (f(k) \leq S_k \leq g(k) : k \leq n) \\[2ex]
& \leq & -{\pi^2 \over 8} 
    (M(n) - 1) L (1+2\ee)^{-2} \\[2ex]
& \leq & \ee + 4(-{\pi^2 \over 8} + \ee) V (1+2\ee)^{-1} \sum_{k=1}^n 
    (f(k) - g(k))^{-2} 
\end{eqnarray*}
for large $n$.  Letting $\ee \rightarrow 0$ proves that the limsup
in~(\ref{eqth1}) is at most $-{\pi^2 \over 8} V$ .

Proving that the liminf is at least $-{\pi^2 \over 8} V$ is almost identical.  
Fixing $L, \ee > 0$, choose $N$ and $m_0$ as before and this time define $R_k$
to lie between $f$ and $g$ instead of covering them:
$$R_k =  \{ (x,y) : m_k \leq x \leq m_{k+1} \mbox{ and } y_k - (1- 2\ee) d_k 
\leq y \leq y_k + (1- 2\ee) d_k \} .$$  
Let $G_k$ be the event that $(i,S_i) \in R_k$ for $m_k \leq i \leq m_{k+1}$
and that $y_k - {1 \over 4} d_k \leq S_{m_{k+1}} \leq
y_k + {1 \over 4} d_k$. Since $m_{M(n) + 1} \geq n$, the probability 
we are trying to bound 
from below is at least $\P (g(k) \leq S_k \leq f(k) : k = 1, \ldots , 
m_{M(n) + 1})$, which may be written as
$$\P (g(k) \leq S_k \leq f(k) : k = 1, \ldots , m_0) \prod_{i=0}^{M(n)}
   \P (G_k \| S_{m_k} , G_{k-1}) .$$
When $G_{k-1}$ occurs, the value of $S_k$ is certainly between 
$y_k - {1 \over 2} d_k$ and $y_k + {1 \over 2} d_k$, so 
$$\P (G_k \| S_{m_k} , G_{k-1}) \geq \inf_{y_k - {1 \over 2} d_k \leq
x \leq y_k + {1 \over 2} d_k} \P (G_k \| S_k = x).$$  
Now rescaling each rectangle,  applying~(\ref{lower rw}) with $c_1 = 
-(1/2)(1-2\ee)^{-2} , c_2 = (1/2) (1-2\ee)^{-2} , c_3 =
-1/4, c_4 = 1/4$, and taking limits establishes that the liminf 
in~(\ref{eqth1})
is at least $-{\pi^2 \over 8} V$, finishing the proof of the lemma.    

Finally, the corollary is proved by letting $f(n) = c_2 n^{1/3}, g(n)
= c_1 n^{1/3}$, and verifying the flatness hypothesis; 
summing $[f(k) - g(k)]^{-2}$ from 1 to $n$ gives $(3+o(1)) 
(c_2 - c_1)^{-2} n^{1/3}$ and the desired conclusion follows.    $\Cox$

The next step is to apply these random walk estimates to prove
a tree-indexed version of Kesten's result on branching Brownian
motion staying above zero.

\begin{lem} \label{lem small prob}
Let $\Gamma, X_\sigma , S_\sigma$ and $\beta$ be as above.
Assume that $X$ is top-heavy and let $\lambda_0$ be the value of 
$\lambda$ minimizing $\E e^{\lambda X}$
(which must exist and be strictly less than one,
according to the definition of the term top-heavy).
There exists  a positive real number $c$, independent of $\Gamma$,
such that if $|\Gamma_n| \leq e^{c n^{1/3} + n \beta}$ for all $n$, then
\begin{equation} \label{U is small}
\P (\max_{\sigma \in \Gamma_n} \min_{\tau \leq \sigma}
   S(\tau) \geq -2 c(1-\lambda_0)^{-1} n^{1/3}) \rightarrow 0
\end{equation}
as $n \rightarrow \infty$.  In other words with high probability,
for each sufficiently large $n$, no path from $\rho$ of length $n$ 
stays above $-2 c (1-\lambda_0)^{-1} n^{1/3}$.
\end{lem}

To establish this, first record some elementary facts about large 
deviations.
\begin{pr} \label{pr large dev}
Let $X_n$ be i.i.d.\ and $S_k = \sum_{i=1}^k X_i$.  Let $\beta (X)$
be the backwards push, and let $\lambda_0 (X)$ be the $\lambda \in 
[0,1]$ minimizing $\E e^{\lambda X}$.  
Then the following three inequalities hold.

\begin{quotation}
$(i)$  ~For any real $u$,
$$ \P (S_n \geq u) \leq e^{-\beta n -\lambda_0 u} ; $$

$(ii)$  ~For any real $y$,
$$ \E e^{S_n} I(S_n \leq y) \leq (1 - \lambda_0)^{-1} e^{(1-\lambda_0) 
   y - \beta n} .  $$

$(iii)$  ~For any real $y$,
$$ \E e^{S_n \wedge y} \leq [1+(1 - \lambda_0)^{-1}] e^{(1-\lambda_0) 
   y - \beta n} .  $$
\end{quotation}
\end{pr}

\noindent{{\bf Proof:}}  The first claim is just Markov's inequality:
$$\P (\S_n \geq u) \leq e^{-\lambda_0 u} \E e^{-\lambda_0 S_n}.$$
For the second claim, integrate the first by parts:
\begin{eqnarray*}
\E (e^{S_n} I(S_n \leq y)) & = & \int_{-\infty}^y e^u \P (S_n \in du) \\[2ex]
& = & \int_{-\infty}^y e^u \P (u \leq S_n \leq y) du \\[2ex]
& \leq & \int_{-\infty}^y e^u \P (u \leq S_n) du \\[2ex]
& \leq & \int_{-\infty}^y e^u e^{-\lambda_0 u - \beta n} du \\[2ex]
& = & (1-\lambda_0)^{-1} e^{(1 - \lambda_0)y - \beta n} .
\end{eqnarray*}
Finally, the third claim follows from the first two, using 
$$\E e^{S_n \wedge y} = \E e^{S_n} I(S_n \leq y) + e^y \P (S_n > y) .$$
It should be remarked that more careful estimates give an extra factor 
of $(1+o(1))(2 \pi n \E Y^2)^{-1/2}$ in the RHS of each inequality
which is then asymptotically sharp.    $\Cox$

Next, plug this into a first moment calculation to establish:
\begin{pr} \label{pr unlikely}
Let $\Gamma , X_\sigma , S_\sigma$ and $\beta$ be as above.
Suppose that $|\Gamma_n| \leq  e^{cn^{1/3} 
+ n \beta}$ for some $c > 0$.  Then for any $\ee > 0$,
$$\P (S(\sigma) \geq (1+\ee) c \lambda_0^{-1} n^{1/3} \mbox{ for some
$\sigma$ with } |\sigma| \leq n) \rightarrow 0 $$
as $n \rightarrow \infty$.
\end{pr}

\noindent{{\bf Proof:}}  For each fixed $L$ it is clear that 
\begin{equation} \label{eqxx}
\P (S(\sigma) \geq (1+\ee) c \lambda_0^{-1} n^{1/3} \mbox{ for some
$\sigma$ with } |\sigma| \leq L) 
\end{equation}
goes to zero as $n \rightarrow \infty$.  On the other hand,
\begin{eqnarray*}
 && \P (S(\sigma) \geq (1+\ee) c \lambda_0^{-1} n^{1/3} \mbox{ for some
    $\sigma$ with } n \geq |\sigma| > L) \\[2ex] 
& \leq & \P (S(\sigma) \geq (1+\ee) c \lambda_0^{-1} |\sigma|^{1/3} 
    \mbox{ for some $\sigma$ with } n \geq |\sigma| > L) \\[2ex]
& \leq & \sum_{m>L} \P (S(\sigma) \geq (1+\ee) c \lambda_0^{-1} m^{1/3}
    \mbox{ for some } \sigma \in \Gamma_m) .
\end{eqnarray*}
For $\sigma \in \Gamma_m$, Proposition~\ref{pr large dev} part $(i)$ 
implies 
$$\P (S(\sigma) \geq (1+\ee) c \lambda_0^{-1} m^{1/3}) \leq 
   e^{-\lambda_0 (1+\ee) c \lambda_0^{-1} m^{1/3} - \beta m} .$$
Multiplying by $|\Gamma_m|$ gives
$$\P (\S(\sigma) \geq (1+\ee) c \lambda_0^{-1} m^{1/3} 
   \mbox{ for some } \sigma \in \Gamma_m) \leq  e^{- \ee c n^{1/3}}. $$
This is summable in $m$, so the sum over $m > L$ goes 
to zero as $L \rightarrow \infty$, which together with~(\ref{eqxx}) 
proves the proposition.   $\Cox$ 

\noindent{{\bf Proof of Lemma}}~\ref{lem small prob}:  
Let $\mu$ be the common
distribution of the $X(\sigma)$ and let $Y_1 , \ldots , Y_n$ be $\iid$
random variables whose law $\mu'$ satisfies 
\begin{equation} \label{tilted}
{ d \mu' \over d\mu} (x) = e^{\lambda_0 x} / \E e^{\lambda_0 X}
   = e^{\lambda_0 x + \beta (\mu)} .
\end{equation}
Informally, $\mu'$ is $\mu$ tilted in the large deviation sense so
as to have mean zero.  
The assumption that $\mu$ is top-heavy by definition implies that
$\E e^{\lambda X} < \infty$ for $\lambda$ in some neighborhood of
$\lambda_0$, hence $\E e^{\lambda Y_1} < \infty$ for $\lambda$
in some neighborhood of zero, and in particular $V \deq \E Y_1^2 < \infty$.

Choose a positive real $c$ for which 
\begin{equation} \label{c is small}
c + 2 \lambda_0 (1- \lambda_0)^{-1} c -{\pi^2 \over 8} 
    {3V \over (2c\lambda_0^{-1} + 2c(1-\lambda_0)^{-1})^2}  < 0 .
\end{equation}
Let 
$$A_n = \left \{ \max_{\sigma \in \Gamma_n} \min_{\tau \leq \sigma} S(\tau)
   \geq -2c(1-\lambda_0)^{-1} n^{1/3} \right \} $$
be the event in~(\ref{U is small}).  Let $G_n$ be the event 
$$ \bigcup_{|\sigma| \leq n} \{ S(\sigma) \geq 2c 
   \lambda_0^{-1} n^{1/3} \} $$
and let $H_n = A_n \setminus G_n$.
Proposition~\ref{pr unlikely} shows that $\P (G_n) \rightarrow 0$ so
to show that $\P (A_n) \rightarrow 0$ it suffices to show that
$\P (H_n) \rightarrow 0$.  

To see this, fix $\sigma \in \Gamma_n$ and write $\P (H_n) \leq
|\Gamma_n| Q_n$ where 
$$Q_n = \P (2 c \lambda_0^{-1} n^{1/3} \geq S(\tau) \geq -2c
   (1-\lambda_0)^{-1} n^{1/3} \mbox{ for all } \tau \leq \sigma) .$$
Let $\nu$ be the law in $\R^n$ of the sequence $(S(\sigma_1) ,
\ldots , S(\sigma_n))$, where $\sigma_1 , \ldots , \sigma_n$
is the path from $\rho = \sigma_0$ to $\sigma = \sigma_n$; of
course $\nu$ is just the law in $\R_n$ of a random walk whose
steps have law $\mu$.  Recalling the tilted variables $Y_n$,
write $T_n = \sum_{i=1}^n Y_i$.  Let $\nu'$ denote the law in 
$\R_n$ of $(T_1 , \ldots , T_n)$ and observe that
$${ d \nu' \over d\nu} (s_1 , \ldots , s_n) = e^{\lambda_0 s_n} / 
  \E e^{\lambda_0 S_n} = e^{\lambda_0 s_n + n \beta} .$$
Use this to get an upper bound on $Q_n$ as follows.
\begin{eqnarray*}
&& Q_n \\[2ex]
& = & \int I(2 c \lambda_0^{-1} n^{1/3} \geq s_i \geq -2c
   (1-\lambda_0)^{-1} n^{1/3} \mbox{ for all } i \leq n) \\
&&  d\nu (s_1 , \ldots , s_n) \\[2ex] 
& = & \int I(2 c \lambda_0^{-1} n^{1/3} \geq s_i \geq -2c
   (1-\lambda_0)^{-1} n^{1/3} \mbox{ for all } i \leq n) \\
&&  \hspace{.5in}{ d\nu \over d\nu'} (s_1 , \ldots , s_n)   
   d\nu' (s_1 , \ldots , s_n) \\[2ex]
& \leq & \int I(2 c \lambda_0^{-1} n^{1/3} \geq s_i \geq -2c
   (1-\lambda_0)^{-1} n^{1/3} \mbox{ for all } i \leq n) \\[1ex]
&& ~~~~~~~~~~\left [ \sup_{s_n 
   \geq -2c(1-\lambda_0)^{-1}n^{1/3}} { d\nu \over d\nu'}
   (s_1 , \ldots , s_n) \right ]  \; d\nu (s_1 , \ldots , s_n) \\[2ex]
& = & \int I(2 c \lambda_0^{-1} n^{1/3} \geq s_i \geq -2c
   (1-\lambda_0)^{-1} n^{1/3} \mbox{ for all } i \leq n) \\[1ex]
&& \hspace{1in}  \exp (2c\lambda_0 
    (1-\lambda_0)^{-1} n^{1/3} - n\beta) \; d\nu (s_1 , \ldots , s_n)
    \\[2ex]
& \leq &  \exp (2c\lambda_0 (1-\lambda_0)^{-1} n^{1/3} - n\beta) \\ 
&& \hspace{.5in} \P (2c\lambda_0^{-1} n^{1/3} \geq T_i \geq -
   2c(1-\lambda_0)^{-1} \mbox{ for all } i \leq n) \\[2ex]
& \leq & \exp \left (- n\beta + n^{1/3} [{2c\lambda_0 \over 1-\lambda_0}
   - {\pi^2 \over 8} {3V \over (2c\lambda_0^{-1} + 
   2c(1-\lambda_0)^{-1})^2} + o(1) ] \right )
\end{eqnarray*}
by Corollary~\ref{cor1}.  Thus $\P (H_n) \leq |\Gamma_n| Q_n \leq
\exp (cn^{1/3} + n \beta) Q_n \leq$
$$\exp \left ( n^{1/3} \left [ c  + 2c\lambda_0 (1-\lambda_0)^{-1}
    - {\pi^2 \over 8} {3V \over (2c\lambda_0^{-1} + 2c(1-\lambda_0)^{-1})^2}
    + o(1) \right ] \right ) .$$
By choice of $c$, this is $\exp ((K + o(1)) n^{1/3})$ for some $K<0$, 
so $\P (H_n) \rightarrow 0$, proving the lemma.   $\Cox$

\section{Proof of the main theorem}

The case where $X \equiv -\beta$ is done in \cite{Ly1}.  For
part~$(i)$ of the nondegenerate case, use Lemma~\ref{bottleneck}, 
showing that
$\sum_{|\sigma| = n} U(\sigma)$ goes to zero in probability 
by computing a truncated expectation.
Let $c_1$ be the constant $c$ from Lemma~\ref{lem small prob}
and let $U_n = \sum_{|\sigma| = n} U(\sigma)$.  Let $G_n$ be the 
event that $\max_{|\sigma| = n} U(\sigma) \geq \exp (-2c_1
(1-\lambda_0)^{-1} n^{1/3})$.  Then
\begin{eqnarray*}
P (U_n > \ee) & \leq & \P (G_n) + \P (U_n > \ee 
   \mbox{ and } G_n^c) \\[2ex]
& \leq & \P (G_n) + \ee^{-1} \E U_n I(G_n^c).
\end{eqnarray*}
Lemma~\ref{lem small prob} showed that $\P (G_n) \rightarrow 0$, so it 
remains to show that for any $\ee$, $\E (U_n I(G_n^c)) \rightarrow 0$. 
 
Observe that $U_\sigma I(G_n^c) \leq \exp (S_\sigma \wedge 
-2c (1-\lambda_0)^{-1} n^{1/3})$.  Hence for $\sigma \in \Gamma_n$,
\begin{eqnarray*}
&& \E U_n I(G_n^c) \\[2ex]
& = & |\Gamma_n| \E U_\sigma I(G_n^c) \\[2ex]
& \leq & |\Gamma_n| \E \exp (S_\sigma \wedge -2c_1 (1-\lambda_0)^{-1}
   n^{1/3} ) \\[2ex]
& \leq & |\Gamma| [1+(1-\lambda_0)^{-1}] \exp ((1- \lambda_0)
   (-2c_1(1-\lambda_0)^{-1} n^{1/3}) - n \beta) \\[2ex]
&& \mbox{ by Proposition~\protect{\ref{pr large dev}}, part } (iii), \\[2ex]
& \leq &  [1+(1-\lambda_0)^{-1}] \exp (-c_1 n^{1/3})
\end{eqnarray*}
by the assumption on $|\Gamma_n|$.  This goes to zero, thus
$U_n \rightarrow 0$ in probability, proving 
part $(i)$ of Theorem~\ref{th2.3}.     

Part~$(ii)$ is proved by exhibiting an infinite path
along which the resistances $e^{-S_\sigma}$ are 
summable.  In fact the proof finds an infinite path along which
$S(\sigma) / c |\sigma|^{1/3}$ is bounded above and below.

Pick any $c, L, \ee > 0$ and any $c_2 = K + 2c \lambda_0 + ({\pi^2 \over 8}) 
{3V \over (.9 c)^2}$, where $M$ shall be chosen later.  Define $W$ to be 
the random subset of vertices $\sigma$ of $\Gamma$ with the property that
for every $\tau \leq \sigma$ with $|\tau| > L$, 
$${c|\tau|^{1/3} \over 10} \leq S_\tau \leq c |\tau|^{1/3} ,$$
where $L$ is large enough so that 
$W$ intersects each $\Gamma_n$ with positive probability.
To prove the theorem, it  suffices to show that $W$ is infinite with positive 
probability; this follows from Lemma~\ref{lem2mm} and the hypothesis 
of the theorem, provided that
\begin{equation} \label{QB2B}
\sup_n a(n,k)/a(n,n)^2 \leq e^{c_2 k^{1/3} + k \beta}
\end{equation}
for all but finitely many $k$, where $a(n,k) = \P (\rho \leftrightarrow
\sigma , \tau)$ for vertices $\sigma , \tau \in \Gamma_n$ with
$|\sigma \wedge \tau| = k$.  

To establish~(\ref{QB2B}), begin with
\begin{equation} \label{QB3B}
a(n,k) / a(n,n)^2 = \P (\sigma \wedge \tau \in W)^{-1}\, {\P (\sigma \in W \|
   \tau \in W) \over \P (\sigma \in W \| \sigma \wedge \tau \in W)} .
\end{equation}
Recall the tilted random variables $Y_n$ and $T_n$ whose law $\mu'$ is 
defined by~(\ref{tilted}).  

Fix any $k \leq n$ and $\sigma , \tau \in \Gamma_n$ with $|\sigma \wedge
\tau| = k$.  Let $C(a,b)$ denote the set of sequences $(s_a , \ldots , s_b) 
\in \R^{b-a+1}$ for which $c s_j^{1/3} / 10 \leq s_j \leq c s_j^{1/3}$ for all $j \in [a,b]$.  Let $\nu$ and $\nu'$ respectively denote
the law of $(S_1 , \ldots , S_k)$ and $(T_1 , \ldots , T_k)$ and
let $\nu_y$ and $\nu_y'$ denote the laws of $(S_k , \ldots , S_n)$
and $(T_k , \ldots , T_n)$ conditioned respectively on $S_k = y$
and $T_k = y$.  Write the first factor on the RHS of~(\ref{QB3B}) as
$$ \left [ \int I (\alpha \in C(1,k)) \, d\nu (\alpha) \right ]^{-1} .$$
Changing the integrating measure to $\nu'$ yields
\begin{eqnarray}
&& \left [ \int I (\alpha \in C(1,k)) \, d\nu' (\alpha) \, { d\nu(\alpha) \over
   d\nu' (\alpha)} \right ]^{-1} \nonumber \\[2ex]
& \leq & \sup_{\alpha \in C(1,k)}  { d\nu' (\alpha) \over d\nu (\alpha)} 
   \left [ \int I (\alpha \in C(1,k)) \, d\nu' (\alpha) \right ]^{-1} 
   \nonumber \\[2ex]
& = & e^{\lambda_0 c k^{1/3} + k \beta} (\nu' (C(1,k)))^{-1}  
   \nonumber \\[2ex]
& = & \exp \left [ \lambda_0 c k^{1/3} + k \beta + k^{1/3}  
   ({\pi^2 \over 8} {3V \over (.9 c)^2} +o(1)) \right ]  \label{notstupid}
\end{eqnarray}
by Corollary~\ref{cor1}, since $\nu'$ is the law of a mean zero, finite 
variance random walk, and it has been assumed that $\nu' (C(1,k))$ 
never vanishes.  

For the second factor on the RHS of~(\ref{QB3B}),
let $\mu_1$ be the law of $S(\sigma \wedge \tau)$
conditional on $\tau \in W$ and let $\mu_2$ be the law of 
$S(\sigma \wedge \tau)$
conditional on $\sigma \wedge \tau \in W$.  The second factor is then
\begin{equation} \label{second factor}
{\int \left [ \int I (\alpha \in C(k,n)) \, d\nu_y (\alpha)
     \right ] \, d\mu_1 (y)    \over
 \int \left [ \int I (\alpha \in C(k,n)) \, d\nu_y (\alpha)
     \right ] \, d\mu_2 (y) } \; .
\end{equation}
Changing the integrating measure again, this becomes
\begin{eqnarray*}
&& {\int \int I(\alpha \in C(k,n)) \, {d\nu \over d\nu'} (\alpha)
    \, d\nu_y' (\alpha)  d\mu_1 (y) \over
    \int \int I(\alpha \in C(k,n)) \, {d\nu \over d\nu'} (\alpha)
    \, d\nu_y' (\alpha)  d\mu_2 (y)} \\[2ex]
& \leq & {\int  I(\alpha \in C(k,n)) d\nu_y'(\alpha) d\mu_1 (y) \over
          \int  I(\alpha \in C(k,n)) d\nu_y'(\alpha) d\mu_2 (y)} \;
    {\sup_{\alpha \in C(k,n)} {d\nu \over d\nu'} (\alpha) \over 
     \inf_{\alpha \in C(k,n)} {d\nu \over d\nu'}  (\alpha)} \\[3ex]
& \leq & {\sup_y \nu_y' (C(k,n)) \over \inf_y \nu_y' (C(k,n))} \;
    {\sup_{\alpha \in C(k,n)} {d\nu \over d\nu'} (\alpha) \over 
     \inf_{\alpha \in C(k,n)} {d\nu \over d\nu'}  (\alpha)} \\[3ex]
& = & {\sup_y \nu_y' (C(k,n)) \over \inf_y \nu_y' (C(k,n))} \;
    e^{\lambda_0 c k^{1/3}} .
\end{eqnarray*}

The argument is then finished by establishing 
\begin{equation} \label{stupid}
{\sup_y \nu_y' (C(k,n)) \over \inf_y \nu_y' (C(k,n))} \leq e^{M k^{1/3}}
\end{equation}
for some $M > 0$, since then multiplying inequalities~(\ref{notstupid})
and~(\ref{stupid}) bounds $a(n,k)/a(n,n)^2$ from above by 
$\exp \left [ 2\lambda_0 c k^{1/3} + k \beta + k^{1/3} ({\pi^2 \over 8} 
{3V \over (.9 c)^2} + M) \right ]$, which is at most $e^{c_2 k^{1/3} 
+ k \beta}$ by choice of $c_2$, yielding~(\ref{QB2B}).

It remains to establish~(\ref{stupid}).  An argument is given
for the case where the distribution of the $X(\sigma)$'s has
an absolutely continuous component near zero, the lattice
case being similar.  By hypothesis, the measure $\mu'$ 
has density at most $A$ and is greater than some constant, $a > 0$,
times Lebesgue measure on some interval $(-b,b)$.  Call this latter
measure $\pi$.  Let $l = k + c k^{1/3} / b$.  Write 
$$\nu_y' (C(k,n)) = \int m_z (C(k+l , n)) \, d m^y (z)$$
where $m_z$ is the law of $(T_{k+l} , \ldots , T_n)$ conditioned
on $T_{k+l} = z$ and $m^y$ is the (deficient) law of $T_{k+l}$
conditioned on $T_k = y$ and killed if $T_i \notin
[c i^{1/3} / 10 , c i^{1/3}]$ for some $k \leq i \leq k+l$. 
It is easy to see that 
$$m^y \geq \pi^y \geq a^l \mbox{ times Lebesgue measure on }
   [c (k+l)^{1/3} / 10 , c (k+l)^{1/3}]$$
where $\pi^y$ is the measure on sequences $s_k , \ldots , s_{k+l}$
with $s_k = y$, having increments distributed as $\pi$ and killed 
if $s_i \notin [c i^{1/3} / 10 , c i^{1/3}]$ for some $k \leq i 
\leq k+l$. Then for any $x,y \in [ck^{1/3} / 10 , ck^{1/3}]$, 
\begin{eqnarray*}
&& \nu_y' (C(k,n)) \\[2ex]
& = & \int m_z (C(k+l , n)) \, d m^y (z) \\[2ex]
& \geq &  a^l \int m_z (C(k+l , n)) \, d \lambda (z) \\[2ex]
& \geq &  (a/A)^l \int m_z (C(k+l , n)) \, d m_x (z) \\[2ex]
& = &  (a/A)^l \nu_x' (C(k,n)) . 
\end{eqnarray*}
Checking against the value of $l$ proves~(\ref{stupid}) and
the theorem.   $\Cox$

\renewcommand{\baselinestretch}{1.0}\large

\noindent{Keywords:} RWRE, critical RWRE, tree, tree-indexed random walk

\noindent{Subject classification: } Primary: 60J15.  Secondary: 60G60, 60G70, 
60E07.


\normalsize
\begin{thebibliography}{YMN}

\bibitem{BP}
Benjamini, I. and Peres, Y. (1993).  Markov chains indexed by a tree.
{\em Ann. Probab. to appear.}

\bibitem{DS}
Doyle, P. and Snell, J. L. (1984).  Random walks and electrical networks.
Mathematical Association of America: Washington.

\bibitem{Du}
Durrett, R. (1984).  Brownian motion and martingales in analysis.
Wadsworth: Monterey, CA.

\bibitem{Ke}
Kesten, H. (1978).  Branching Brownian motion with absorption.
{\em Stoch. Pro. Appl.} {\bf 7} 9 - 47.

\bibitem{Ly1}
Lyons, R. (1990).  Random walks and percolation on trees.  {\em Ann.
Probab.} {\bf 18} 931 - 958.

\bibitem{Ly2}
Lyons, R. (1992).  Random walks, capacity and percolation on trees.
{\em Ann. Probab. to appear.}

\bibitem{LP}
Lyons, R. and Pemantle, R. (1992).  Random walk in a random environment and
first-passage percolation on trees.  {\em Ann. Probab.} {\bf 20}  125 - 136.

\bibitem{Pe}
Pemantle, R. (1988).  Phase transition in Reinforced random walk and 
RWRE on trees.  {\em Ann. Probab.} {\bf 16} 1229 - 1241.

\bibitem{PP}
Pemantle, R. and Peres, Y. (1991).  Critical random walk in random
environment on a tree.  {\em Preprint.}

\end{thebibliography}
\end{document}